\documentclass[12pt,a4paper,twoside]{article}

\usepackage{amsmath,amsthm,amsfonts}

\usepackage{fancyheadings}

\newif\ifpdf
\ifx\pdfoutput\undefined
\pdffalse 
\else
\pdfoutput=1 
\pdftrue
\fi

\ifpdf
\usepackage[pdftex]{graphicx}
\else
\usepackage{graphicx}
\fi

\makeatletter
\renewcommand{\ps@plain}{%
\renewcommand{\@oddhead}{}%
\renewcommand{\@evenhead}{}%
\renewcommand{\@evenfoot}{}%
\renewcommand{\@oddfoot}{}%
}
\makeatother

\def\back{\noindent\kern-5mm} 
\def \Liminf{\mathop{\underline{\lim}}\limits}
\def\de{{\rm d}}

\def \Label{\label}
\def \bE{{\mathbf E}}

\def \bR{{\mathbb R}}
\def \de{{\rm d}}

\def \mS{{\mathcal S}}
\def \mP{{\mathcal P}}
\def \mR{{\mathcal R}}
\def \mE{{\mathcal E}}
\def \mM{{\mathcal M}}

\def \ra{{\rightarrow}}

\def \lgraf{\{}
\def \rgraf{\}}
\def \zs#1{_{\lower 3pt \hbox{$\scriptstyle#1$}}}

\def\beq{\begin{equation}}
\def\endeq{\end{equation}}
\def\beqn{\begin{eqnarray*}}
\def\endeqn{\end{eqnarray*}}

\newtheorem{theorem}{Theorem}

\newtheorem{definition}{Definition}

\newtheorem{example}{Example}

\def\thetitle{Efficiency of a class of unbiased estimators  for the invariant distribution function of a diffusion process \footnote{This work has been partially supported by the local Ggrant sponsored by University of Bergamo: {\it Theoretical and computational problems in statistics for
continuously and discretely observed diffusion processes} and MIUR 2004 Grant. } }

\title\thetitle
\author{
Ilia Negri\footnote{%
Department of Management and Information Technology,
Viale Marconi 5,
24044 Dalmine (BG), Italy. 
{\tt ilia.negri@unibg.it}}
}

\begin{document}

\date{}

\ifpdf
\DeclareGraphicsExtensions{.pdf, .jpg}
\else
\DeclareGraphicsExtensions{.eps, .jpg}
\fi

\maketitle

\begin{abstract}
We consider the problem of the estimation of the invariant distribution function of an ergodic diffusion process when the drift coefficient is unknown.  The empirical distribution function is a natural estimator which is unbiased, uniformly consistent and efficient in different metrics. Here we study the properties of optimality for an other kind of estimator recently proposed. We consider a class of unbiased estimators  and we show that they are also  efficient in the sense that their asymptotic risk, defined as  the integrated mean square error, attains an asymptotic minimax lower bound. 
\end{abstract}

\par

\par

\vfill
\eject

\pagestyle{fancy}

\baselineskip 22pt

\section{Introduction}
We consider the problem of the estimation of the distribution function
$F(x)$, $x \in \bR$
by the observation of a diffusion process $\lgraf X_t : 0\leq t \leq T
\rgraf $.
We suppose that the process $X_t$, $t \geq 0$ possesses the ergodic property
with invariant measure $\mu$ and $F(x)=\mu \lgraf (-\infty, x] \rgraf$.
A natural estimator for $F(x)$,  $x \in \mathbb R$,   is the empirical distribution
function
$\hat{F}_T(x)$.
It is well known that this estimator is uniformly
consistent by the Glivenko--Can\-tel\-li theorem and asymptotically normal (Kutoyants, 1997).
The problem of the asymptotically efficiency of the empirical distribution function has been considered for different model and different metrics. For the model of independent and identical distributed random variables the empirical distribution function is asymptotically efficient, in a global framework, in the sense that  its integrated mean square error attains the lower bound given for all the estimators of the  distribution function.  Such result has been established earlier by Levit (1978) and Millar
(1979) using the theory of local asymptotic normality. Gill and
Levit (1995) obtained the same result using a different approach
based on a multidimensional version of the van Trees inequality.
The same approach introduced by Gill and Levit was successfully  applied in Kutoyants and Negri, (2001) to prove that the empirical distribution function is asymptotically efficient in the problem of invariant distribution estimation for ergodic diffusion processes. For the same model Negri (1998) has proved the asymptotically efficiency of the empirical distribution function when the 
metric utilized in the risk function is based on the {\em sup} norm. 

Recently (Kutoyants, 2004) a class of unbiased estimator for the invariant distribution function has been introduced.  These estimators, that do not contain the empirical distribution function as particular case,  are consistent and asymptotically normal. In this work we prove that they are also asymptotically efficient in the sense that their integrated mean square error attain the lower bound given for all the estimators of the invariant distribution function. 
As in the case of the estimation of the invariant density, we have many efficient estimators, so the problem of finding the second order efficient estimator arise naturally (see Dalalyan and Kutoyants 2004 where the problem is considered for the invariant density estimation). This problem it is not considered here, but it will be an argument of future researches. 

The note is organized as follow. In the next section we present the statement of the problem and the
assumptions.
In Section 3 we present the lower bound for the risk. In Section 4 
we prove that the class of unbiased estimator attain this bound and finally,  in Section 5, we give  same examples of such estimators. 

\section{Preliminaries}
\Label{sec2}
In this section we introduce the model and its first properties, while the statistical problem will be presented in the next section. 
Let us consider a one dimensional
diffusion process
\begin{equation}
\de X_t = S(X_t) \de t+\sigma(X_t)\de W_t, \quad
X_0,  \quad  t\geq 0
\Label{eds}
\end{equation}
where $\{W_t:\ t\geq 0\}$ is a standard Wiener process, and the initial value $X_0$ is independent of $W_t$, $t\geq 0$.
The drift coefficient  $S$ will be supposed unknown to the observer
and the diffusion coefficient $\sigma^2$
will be a known positive function. 
Let us introduce the condition:

$\mE\mS$. {\em The function $S$ is locally bounded, the function $\sigma^2$ is positive and continuous and for some $A>0$ the condition
$ xS(x)+\sigma(x)^2\leq A(1+x^2)$, $x\in \bR$ holds.}

Under the condition $\mE\mS$ the equation \eqref{eds} has an unique weak solution (see Durrett, 1996, p. 210). To guarantee ergodicity we introduce the following condition:

$\mR\mP$. {\em The function $S$  and $\sigma$ are such that:
$$
V_S(x)= \int_0^x \exp\left\{-2\int_0^y \frac{S(v)}{\sigma(v)^2}\de v \right\} \de y \to \pm \infty, \quad as \quad x\to \pm \infty
$$
and 
$$
G(S)=\int_{-\infty}^{+\infty} \frac{1}{\sigma(x)^2}
\exp\left\{2\int_0^x \frac{S(v)}{\sigma(v)^2}\de v \right\}\de x < +\infty.
$$
}
If the condition $\mR\mP$ is satisfied then the weak solution of \eqref{eds}, $\{X_t, t\geq 0\}$, has the ergodic property (see for example Gikhman and Skorohod, 1972), that is, there exists an invariant probability measure $\mu_S$ such that for every
measurable function $g$ such
that
$\bE_S |g(\xi)| <\infty$, we have with  probability one,
$$
 \lim_{T\to\infty}\frac1T \int_0^T g(X_t)\de t =
\int_\bR g(z) f_S(z) \de z= \bE_S (g(\xi))
$$
where $\xi$ has the invariant measure as distribution, $\bE_S$ denote the mathematical expectation with respect to $\mu_S$,
and $f_S$ is the invariant density given by
\begin{equation}
f_S(y)=\frac{1}{G(S)\sigma(y)^2}
\exp\left\{2\int_0^y \frac{S(v)}{\sigma(v)^2}\de v \right\}.
\label{invdens}
\end{equation}
Suppose we observe different diffusion processes $\{X_t :  0\leq t\leq T\}$ given by equation \eqref{eds} with drift coefficients respectively given by $S_1$, $S_2$ and $S_0=0$ and initial value respectively $X_0^{1}$, $X_0^{2}$ and $X_0^{0}$. Let us introduce the following condition.

$\mE\mM$. {\em The functions $S_1$, $S_2$  and $\sigma$ satisfy condition $\mE\mS$ and the densities (with respect to the Lebesgue measure)  of  the corresponding initial values $X_0^{1}$, $X_0^{2}$ and $X_0^{0}$ have the same support (if the initial value is nonrandom, then we suppose that it takes the same value for all processes).} 

If  condition $\mE\mM$ holds true,  all the measures $P_S^T$, for different $S$, induced by the process $\{X_t, :  0\leq t\leq T\}$, in the space $\mathcal C_T$, the space of all the continuos function on $[0,T]$ with uniform metric and Borel $\sigma$-algebra $\mathcal B(\mathcal C_T)$, are equivalent. See Kutoyants,  2004, p. 34 and the references therein. 

\section{The asymptotic global bound}
Given the diffusion process \eqref{eds}
we suppose that conditions $\mE \mS$ and $\mR \mP$ are satisfied, that $X_0$ 
has density $f_S$, given by \eqref{invdens}, so the process $\{X_t, t\geq 0\}$ is ergodic and  strictly stationary.  We are interested in the estimation of the invariant distribution function 
\begin{equation}
F_S(x)=\int_{-\infty}^x \frac{1}{G(S)\sigma(y)^2}
\exp\left\{2\int_0^y \frac{S(v)}{\sigma(v)^2}\de v \right\}\de y
\label{sdfS}
\end{equation}
by the observation
$X^T=\{X_t : 0\leq t \leq T\}$ solution of \eqref{eds} when $\sigma$
is known and $S$ is unknown. Let us denote by $\bE_S^T$  the
mathematical expectation with respect to the measure $P_S^T$. 

Let $\bar F_T(x)$ be any estimator of \eqref{sdfS} for $x\in\bR$.
We define the integrated mean square error as
\begin{equation}
\rho_T(\bar F_T, F_S) = T \bE_S^T \int_\bR |\bar F_T(x)-F_S(x)|^2
\nu(\de x)
\end{equation}
where $\nu$ is a finite measure on $\bR$. We will refer  also to it
as {\it global risk}.
A natural estimator of $F_S(x)$ for $x\in \mathbb R$ is the empirical distribution
function defined as follows
\begin{equation}
\hat F_T(x)=\frac1T\int_0^T \chi\zs{\{X_t<x\}} \de t, \quad x\in \mathbb R.
\label{EDF}
\end{equation}
This estimator is uniformly consistent, asymptotically normal and
asymptotically efficient  in
the sense  that the empirical distribution function achieves
a local asymptotic minimax lower bound for the integrated mean square error
of an arbitrary estimator.
For a fixed function $\sigma$ let us introduce the classes
$
{\cal S }_\sigma=\left\{  S : \text{conditions } \mE\mS, \mE \mM, \mR\mP \text{ are fulfilled} \right\}
$ 
and $\cal{S}^*_\sigma \subset \cal{S}_\sigma$ such that for every $S_*$ in $\cal{S}^*_\sigma$ 
there exist a  $\delta>0$, and a vicinity
$V_\delta=\left\{S :
\sup_{x\in\bR}|S_*(x)-S(x)|<\delta, \quad S \in \cal{S}^*_\sigma\right\},
$
such that 
$
\sup_{S\in V_\delta}G(S)<+\infty.
$

For $x$ and $y$ in $\bR$ we denote by $x\wedge y$ and by $x \vee y$
respectively the minimum and the maximum
between $x$ and $y$. Let us introduce the function
$$R_S(x,y)  =
4\int_{-\infty}^{+\infty}
\frac{ F_S(v\wedge x)\big(1-F_S( v \vee x)\big)F_S(v\wedge y)
\big(1-F_S(v \vee y)\big)}{\sigma(v)^2f_S(v)}\de v.
$$
and the quantity
$$\rho_*(S)
 = \int_\bR 4 \bE_S \left(\frac{F_S(x)F_S(\xi)-  F_S(\xi\wedge x)}{\sigma(\xi)f_S(\xi)} \right)^2 \nu(\de x)= \int_\bR R_S(x,x) \nu(\de x)
$$
Let us introduce the following condition.

${\cal Q}_1$. {\em The function $S_*\in {\cal S}_\sigma^*$ and for some $\delta>0$ 
$$\sup_{S\in V_\delta} \rho_*(S)=\sup_{S\in V_\delta}
\int_\bR 4 \bE_S \left(\frac{F_S(x)F_S(\xi)-  F_S(\xi\wedge x)}{\sigma(\xi)f_S(\xi)} \right)^2 \nu(\de x)< +\infty.
$$
}
\noindent
We have the following result (Kutoyants and Negri, 2001).
\begin{theorem}
\Label{Teo1}
Let $S_*\in {\cal S}_\sigma^*$ and condition ${\cal Q}_1$ be fulfilled. Then 
$$\Liminf_{\delta\ra 0}
\Liminf_{T\to\infty}
\inf_{\bar F_T}
\sup_{S\in V_\delta}
\rho_T(\bar F_T, F_S)\geq \rho_*(S^*)
$$
where the $\inf$ is taken over any estimator $\bar F_T$ of $F_S$.
\end{theorem}
\noindent
The definition of asymptotically efficient estimator arise naturally from the above theorem. 

\begin{definition} Let  condition ${\cal Q}_1$ be fulfilled. Then an estimator $\hat F_T$ is called asymptotically efficient if for any $S_*\in {\cal S}_\sigma^*$ we have 
\begin{equation}
\Liminf_{\delta\ra 0}
\Liminf_{T\to\infty}
\sup_{S\in V_\delta}
\rho_T(\hat F_T, F_S) = \rho_*(S^*)
\label{ae}
\end{equation}
\end{definition}
Put
\begin{equation}
H_{x,S}(y)=
2\int_0^y \frac{F_S(v\wedge x)-F_S(v) F_S(x)}{\sigma(v)^2f_S(v)}\de v.
\label{H}
\end{equation}
Let us introduce  the following condition.

${\cal Q}_2$. {\em Suppose that
$$\sup_{S\in V_\delta}
\int_\bR \bE_SH_{x,S}(\xi)^2 \nu(\de x)<+\infty.
$$
}
The following theorem establishes the asymptotic efficiency of the empirical distribution function \eqref{EDF} in the sense of equality \eqref{ae}. 
\begin{theorem} Let conditions ${\cal Q}_1$, ${\cal Q}_2$ hold and $\rho_*(S)$ be continuous in the uniform topology at the point $S^*$, then the empirical distribution function is asymptotically efficient.
\end{theorem}
It is proved in Kutoyants and Negri  (2001).

\section{A class of unbiased estimators}
In this section we consider a class of estimators of $F_S(x)$ recently introduced in Kutoyants (2004) defined, for $x\in \bR$ as
\begin{equation}
\tilde{F}_T(x)=\frac{1}{T} \int_0^T R_x(X_t) \de X_t +\frac{1}{T} \int_0^T N_x(X_t)\de t
\label{ue}
\end{equation}
where $R_x(y)=2\chi_{\{y<x\}} K_x(y)h(y)$,  $N_x(y)=\chi_{\{y<x\}} K_x(y)h'(y)\sigma^2(y)$, 
$K_x(y)=\int_y^x\frac{\de v}{\sigma^2(v)h(v)}$ 
and $h$ is a positive and continuously differentiable function. It can be proved that these estimators, for different functions $R_x$ and $N_x$ are all unbiased, consistent and asymptotically normal for a fixed $x$. Let us suppose the following conditions hold:
\begin{equation}
\bE_S(R_x(\xi)\sigma(\xi))^2<+\infty, \quad \bE_S|N_x(\xi)|<+\infty, \quad \lim_{y\to -\infty} R_x(y)\sigma(y)^2f_S(y)=0
\label{cond1}
\end{equation}
We have the following result  (Kutoyants,  2004).
\begin{theorem}
Let $S \in \cal{S}_\sigma$, $R_S(x,x)<+\infty$, and conditions \eqref{cond1} be fulfilled. Then the estimator $\tilde{F}_T(x)$ is unbiased, consistent and asymptotically normal with variance given by $R_S(x,x)$. 
\label{teo1}
\end{theorem}
Let us define the following function
\begin{equation}
G_{x}(y)= 2 \int_0^y   \chi_{\{ v<x\}}K_x(v)h(v)\de v
\label{G}
\end{equation}
We introduce the following condition.

${\cal Q}_3$. {\em Suppose that
$$\sup_{S\in V_\delta}
\int_\bR \bE_S G_{x}(\xi)^2 \nu(\de x)<+\infty.
$$
}
The following theorem proves that the new class of estimators defined by \eqref{ue} is asymptotically efficient. 
\begin{theorem}
Let conditions ${\cal Q}_1$, ${\cal Q}_2$ and ${\cal Q}_3$ hold true.  Let $\rho_*(S)$ be continuous in the uniform topology at the point $S^*$, then 
\begin{equation}
\Liminf_{\delta\ra 0}
\Liminf_{T\to\infty}
\sup_{S\in V_\delta}
\rho_T(\tilde F_T, F_S) = \rho_*(S^*)
\label{ae2}
\end{equation}
\label{teo2}
\end{theorem}
\begin{proof}
Let us define $c_x(y)=\chi_{\{ y<x\}} K_x(y)[ 2h(y)S(y)+ h'(y)\sigma^2(y)]-F_S(x)$
and $d_x(y)=2\chi_{\{ y<x\}} h(y)K_x(y)\sigma(y)$.
We have the following representation of the empirical process
\begin{equation}
\sqrt{T}\left( \tilde{F}_T(x)-F_S(x) \right)= \frac{1}{\sqrt{T}} \int_0^T c_x(X_t)\de t+ \frac{1}{\sqrt{T}}\int_0^T d_x(X_t) \de W_t.
\label{ep}
\end{equation}
Now we search for a function $M_{x,S}$ such that
\begin{equation}
M_{x,S}'(y)S(y)+\frac12 M_{x,S}''(y) \sigma^2(y)= c_x(y).
\label{eqM}
\end{equation}
Putting $M_{x,S}'=m$, the equation \eqref{eqM} can be rewritten as $m'=\frac{2c}{\sigma^2}- \frac{2S}{\sigma^2} m$ which have solution
\begin{equation}
m(z)=\frac{2}{f(z)\sigma^2(z)} \int_{-\infty}^z c_x(v)f(v)\de v. 
\label{mz}
\end{equation}
Integrating by part the integral in \eqref{mz} and observing that 
$$
\frac{\de}{\de y}\sigma^2(y)f(y)= 2S(y)f(y) 
$$
the function $m$ can be rewritten as 
\begin{equation}
m(z)=2\chi_{\{ z<x\}} h(z)K_x(z) + 2\frac{F(x\wedge z)-F(x)F(z)}{\sigma^2(z)f(z)}.
\label{eqm}
\end{equation}
Choosing the function $M_{x,S}$ such that $M_{x,S}(0)=0$, it has this form
\begin{equation}
M_{x,S}(y)=\int_0^y 2\chi_{\{ z<x\}} h(z)K_x(z) \de z + \int_0^y 2\frac{F(x\!\wedge\! z)-F(x)F(z)}{\sigma^2(z)f(z)}\de z
\label{M1}
\end{equation}
From \eqref{H}, \eqref{G} and \eqref{M1} we can write
\begin{equation}
M_{x,S}(y)= G_{x}(y)+ H_{x,S}(y).
\label{eqM2}
\end{equation}
By the It\^o formula we can write
$$
\de M_{x,S}(X_t)=\left( M_{x,S}'(X_t)S(X_t)+\frac12 M_{x,S}'' (X_t)\sigma^2(X_t)\right)\de t + M_{x,S}'(X_t) \sigma(X_t)\de W_t.
$$
Now we can substitute the Lebesgue integral in \eqref{ep} by means of this formula and the empirical process \eqref{ep} became
$$
\sqrt{T}\! \left(\! \tilde{F}_T(x)\!-\!F_S(x) \! \right)\! \!=\!\!\frac{M_{x,S}(X_T)\!-\! M_{x,S}(X_0)}{\sqrt{T}}+ \frac{1}{\sqrt{T}}\!\! \int_0^T\! \! \! \!2\frac{F_S(x\!\wedge \!X_t)\!-\!F_S(x)F_S(X_t)}{\sigma(X_t)f_S(X_t)} \de W_t.
$$
From \eqref{eqM2} and conditions  ${\cal Q}_2$ and ${\cal Q}_3$ it follows that
\begin{equation}
\lim_{T\to+\infty}\sup_{S\in V_\delta}\int_{\mathbb R}\bE_S\left( \frac{M_{x,S}(X_T)- M_{x,S}(X_0)}{\sqrt{T}}\right)^2 \nu(\de x)=0
\label{final1}
\end{equation}
Moreover by the continuity of $\rho_*(S)$ at the point $S^*$,
 as in Kutoyants and Negri (2001) we can conclude
\begin{equation}
\lim_{\delta\ra 0} \lim_{T\to+\infty}\sup_{S\in V_\delta}\int_{\mathbb R}\bE_S \left(\frac{1}{\sqrt{T}} \int_0^T\!2\frac{F_S(x\wedge X_t)\!-\!F_S(x)F_S(X_t)}{\sigma(X_t)f_S(X_t)} \de W_t \right)^2 \nu(\de x)=\rho_*(S^*)
\label{final2}
\end{equation}
So by \eqref{final1} and \eqref{final2} it follows \eqref{ae2}, and the proof is concluded. 
\end{proof}
\section{Examples}
In Theorem \ref{teo1} and Theorem \ref{teo2}, conditions ${\cal Q}_1$ and ${\cal Q}_2$ involve the function $S$, indeed the model. Example of function $S$ for which conditions ${\cal Q}_1$ and ${\cal Q}_2$ are satisfied can be found in Kutoyants and Negri (2001) and in Kutoyants (2004). Conditions $\eqref{cond1}$ in Theorem \ref{teo1} and ${\cal Q}_3$ in Theorem \ref{teo2} are on the class of estimators.  The class of unbiased estimators defined by \eqref{ue} is very general and is not empty. In this section we will show that for a very large choice of functions $h(\cdot)$ the related estimator $\tilde F$ satisfy condition ${\cal Q}_3$ and conditions \eqref{cond1}. In all this section let us consider for simplicity $\sigma(u)=1$.
\begin{example}
\label{ex2} 
{\rm
Let us consider $h(u)=1+u^{2p}$, $p\geq 1$. We have 
\begin{equation}
K_x(y)=\int_v^x \frac{\de u}{1+u^{2p}} = \tilde{K}_p(x) - \tilde{K}_p(y) 
\label{K1}
\end{equation}
where $ \tilde{K}_p(\cdot)$ is a primitive of $1/h(\cdot)$. We have that for every 
$p\geq 1$, $| K_x(y)|\leq \pi$. 
From \eqref{K1}  we have $|R_x(y)|=| 2\chi\zs{\{y<x\}} K_x(y) h(y)| \leq 2\pi |h(y)|$
and $
\bE_S(R_x(\xi))^2 \leq 2\pi \int_{\mathbb R} (1+y^{2p}) f_S(y) \de y <+\infty
$
because the invariant distribution admits the moment of every order. 
Moreover $N_x(y)= \chi\zs{\{y<x\}} K_x(y) 2p y^{2p-1}$ and we have
$
\bE_S (|N_x(\xi)|)\leq +\infty
$. Finally 
$
\lim_{y\to -\infty} R_x(y) f_S(y)=0 
$ and conditions \eqref{cond1} are satisfied. 
Let now consider 
$
G_x(y) =\int_0^y R_x(v)\de v $.
We have
\begin{equation}
 \bE_s(G_x(\xi)^2) =\int_{-\infty}^{+\infty} \left( \int_0^y R_x(v)\de v \right)^2 f_S(y)\de y <C
\label{G1}
\end{equation}
where $C$ is a constant that does not depend on $S\in V_\delta$. Equation \eqref{G1} imply  
that condition  ${\cal Q}_3$ is satisfied without any further condition on measure $\nu$. 
Note that for $p=1$, the \eqref{K1} became
$
K_x(v)=\int_v^x \frac{\de u}{1+u^2}= \textrm{artg} x -\textrm{artg} v \leq \pi
$ for every $v$ and $x$. 
}
\end{example}
\begin{example}
\label{ex3} 
{\rm
Let us consider $h(u)=e^{\delta u}$, $\delta>0$. We have
$$
K_x(y)=\int_y^x \frac{\de u}{e^{\delta u}}=\frac{e^{-\delta y}-e^{-\delta x} }{\delta} 
$$
and
\begin{equation}
R_x(y)= -2 \chi \zs{\{y<x\}} \frac{e^{\delta( y-x)} }{\delta}.
\label{R2}
\end{equation}
From \eqref{R2} it follows that  $R_x(y)$ is bounded with respect to both the variables $x$ and $y$. 
In virtue of this fact  conditions  \eqref{cond1} and ${\cal Q}_3$ can easily be check. Also in this case we have not  to set any further condition on measure $\nu$. 
}
\end{example}
\begin{example}
\label{ex4} 
{\rm
Let us consider $h(u)=c$, where $c$ is a real constant. We have
$K_x(v)=\frac{x-v}{c}$, $R_x(v)=2 \chi\zs{v<x}(x-v)$ and $N_x(v)=0$. 
In this case the first of the conditions \eqref{cond1} is not satisfied. In fact it follows that
$
\bE_S(R_x(\xi))^2= 4x^2 F_S(x)+C
$, 
where $C$ is a constant depending on the moment values of the invariant distribution. For each fixed $x$, $\bE_S(R_x(\xi))^2$ is finite but when  $x\to +\infty$, also  $\bE_S(R_x(\xi))^2$  goes to infinity. In any case we have
$
G_x(y)= 2x(y\wedge x)-(y\wedge x)^2
$
and $\int_\mathbb R \bE_S(G_x(\xi)^2)\nu(\de x)$ is finite if
the measure $\nu$ admits the moment of order fourth. 
}
\end{example}

We observe that the empirical distribution function does not belong to this class of estimators. Indeed it would be necessary that 
$R_x(y)=0$  and $K_x(y) h(y)=1$ for every $y$ and $x$ belonging to $\mR$. But if such a function $h$ exists it has to depend on $x$ and this is not possible. 

\subsection*{Acknowledgments}
My warmest thanks go to Yury A. Kutoyants  for his useful advices and
clarifying discussions.

\section*{References}
{\leftskip5mm

\back Dalalyan, A.S.,  Kutoyants, Y.A. (2004),  On second order minimax estimation of invariant density for ergodic diffusion. {\em Statist. Decisions}, {\bf 22},  1, 17--41. 

\back Durrett, R.,  (1996),  {\em Stochastic Calculus. A Pratical Introduction}, (CRC Press, New York). 

\back  Gikhman, I.I., Skorohod, A.V. (1972), {\sl Stochastic
        Differential Equations}, (Springer--Verlag, New York).

\back  Gill, R.D., Levit, B.Y. (1995), Applications of the van Trees
inequality: a Bayesian Cram\'er-Rao bound, {\sl Bernoulli}, {\bf 1}+{\bf 2},
59-79.

\back Kutoyants, Y.A. (1997), Efficiency of the empirical
distribution for ergodic diffusion, {\sl Bernoulli}, {\bf 3}, 4,
445-456.

\back Kutoyants, Y.A. (2004), {\sl Statistical Inference for Ergodic Diffusion Processes}, (Springer, New York).

\back  Kutoyants Y.A.,  Negri I. (2001), On $L_2$-efficiency of empiric distribution for
diffusion process, {\it Theory of  Probability and its Applications}, {\bf 46}, 1, 164-169.

\back  Levit, B.Y. (1978), Infinite-dimensional information
inequalities, {\sl Theory Prob. Appl.}, {\bf 23}, 371-377.

\back Millar, P.W. (1979), Asymptotic minimax theorems for the
sample distribution function, {\sl Z. Wahrscheinlichkeitstheorie
Verw. Geb.}, {\bf 48}, 233-252.

\back Negri, I. (1998), Stationary distribution
function estimation for ergodic diffusion process, 
        {\sl Statistical Inference for Stochastic Processes}, {\bf 1}, 61--84.

}

\end{document}